\documentclass[12pt, reqno]{amsart}

\usepackage{amsmath,enumitem}
\usepackage[english,activeacute]{babel}
\usepackage{mathrsfs}
\usepackage[latin1]{inputenc}
\usepackage{amssymb}
\usepackage{amsthm}
\usepackage{graphicx}
\usepackage{array}
\usepackage{pdflscape}
\usepackage{a4wide}
\usepackage{url}
\usepackage{float}
\usepackage{hyperref}
\usepackage{enumerate}
\usepackage{multirow}
\usepackage[capitalize,nameinlink,noabbrev,nosort]{cleveref}
\usepackage{txfonts}
\usepackage[dvipsnames]{xcolor}

\hypersetup{
    pdftoolbar=true,        
    pdfmenubar=true,        
    pdffitwindow=false,     
    pdfstartview={FitH},    
    pdftitle={},
    pdfauthor={},
    pdfsubject={},
    pdfkeywords={},
    pdfnewwindow=true,      
    colorlinks=true,       
    linkcolor=brown,          
    citecolor=brown,        
    filecolor=brown,      
    urlcolor=brown,           
    unicode=true          
}

\numberwithin{equation}{section}

\theoremstyle{plain}
\newtheorem{theorem}{Theorem}[section]
\newtheorem{corollary}[theorem]{Corollary}
\newtheorem{lemma}[theorem]{Lemma}

\theoremstyle{definition}
\newtheorem{definition}[theorem]{Definition}

\title
{Approximation of discontinuous functions by positive linear operators. A probabilistic approach}

\date{\today}

\author{Jos\'{e} A.~Adell}
\address{Departamento de M\'{e}todos Estad\'{i}sticos, Facultad de Ciencias, Universidad de Zaragoza, Spain}
\email{adell@unizar.es}

\author{P. Garrancho}
\address{Departamento de Matem\'aticas, Universidad de Ja\'en, Spain}
\email{pgarran@ujaen.es}

\author{F.J.  Mart\'inez-S\'anchez}
\address{Departamento de Matem\'aticas, Universidad de Ja\'an, Spain}
\email{fjsanche@ujaen.es}

\begin{document}

\begin{abstract}
	We obtain approximation results for general positive linear operators satisfying mild conditions, when acting on discontinuous functions and absolutely continuous functions having discontinuous derivatives. The upper bounds, given in terms of a local first modulus of continuity, are  best possible, in the sense that we can construct particular sequences of operators attaining them. When applied to functions of bounded variation or absolutely continuous functions having derivatives of bounded variation, these upper bounds are better and simpler to compute than the usual total variation bounds. The particular case of the Bernstein polynomials is thoroughly discussed. We use a probabilistic approach based on representations of such operators in terms of expectations of random variables.

\end{abstract}
	\maketitle

\section{Introduction}
\label{introduction}
Denote by $\mathbb{N}$ the set of positive integers and by $I$ a real interval. It is known (c.f. \cite{adellcal95} and \cite{adellcal96}) that many sequences $L_n$ of positive linear operators  allow for a probabilistic representation of the form
\begin{equation}\label{Ln}
L_n(f,x) = \mathbb{E} \, f \left( x + \frac{Z_n(x)}{\sqrt{n}} \right), \quad x\in I,\  n\in\mathbb{N},
\end{equation}
where $\mathbb{E}$ stands for  mathematical expectation, $Z_n(x)$ is a random  variable   such that $x + Z_n (x) / \sqrt{n} $ takes values in $I$, and $f: I \to\mathbb{R}$ is any measurable function for which the   expectations in  (\ref{Ln}) are finite.

Consider, for instance, the Bernstein polynomials $B_n$ defined by
\begin{equation}\label{2}
 B_n(f,x)=\sum_{k=0}^{n}f\left(\frac{k}{n}\right)\binom{n}{k}x^k(1-x)^{n-k},\quad x\in  [0,1],\ n\in\mathbb{N}.
\end{equation}
Let $ x \in[0,1]$ and $n\in\mathbb{N}.$ To give a probabilistic representation for $B_n$, let $(U_k)_{k\geq 1}$ be a sequence of independent copies of a random variable $U$ having the uniform distribution on $[0,1]$. Define

\begin{equation}\label{2*}S_n (x) = \sum_{k=1 }^n 1_{[0,x]} (U_k),\end{equation}
where $1_A$ denotes the indicator function of the set $A$. Clearly, $S_n (x)$ has
the  binomial law  with parameters $n$ and $x$, that is,
\begin{equation}\label{leybinomial}
  P(S_n(x)=k)=\binom{n}{k}x^k(1-x)^{n-k},\quad k=0,1,\ldots, n.
\end{equation}
We therefore have from (\ref{2})
\begin{equation}\label{BnE}
  B_n(f,x)=\mathbb{E} f \left( \frac{S_n (x)}{n}\right)=\mathbb{E} \; f \left(  x + \frac{Z_n (x)}{\sqrt{n}} \right),
\end{equation}
where
\begin{equation}\label{5}
  Z_n(x)=\frac{S_n(x)-nx}{\sqrt{n}}.
 \end{equation}

Fix $x\in \overset{\circ}{I}$, the interior set of $I$. Let $B(I)$ be the set of all measurable bounded real functions defined on $I$ whose right and left limits at $x$, respectively denoted by $f(x+)$ and $f(x-)$, are finite. Also, denote by $DB(I)$ the set of absolutely continuous functions $\phi$ such that
\begin{equation}\label{phi}
\phi (y) = \phi (x) + \int_x ^y f(u)\, du,\quad y\in I, \ f\in B(I).
\end{equation}

The aim of this paper is to give estimates for both
\begin{equation*}
L_n(f,x)-\frac{f(x+)+f(x-)}{2} \mbox{\qquad and\qquad } L_n(\phi,x)-\phi(x),
\end{equation*}
where $L_n$ is a positive linear operator of the form (\ref{Ln}) satisfying the following two properties
\begin{equation}\label{Prop1}
P(|Z_n(x)| \geq \theta ) \leq C_n (\beta) \exp(-\beta \theta),\quad \theta \geq 0,
 \end{equation}
 for some positive constants $\beta$ and $C_n(\beta)$, possibly depending upon $x$, and
\begin{equation}\label{Prop2} \lim_{n\rightarrow\infty}P(Z_n(x)>0)=\lim_{n\rightarrow\infty}P(Z_n(x)<0)=\frac{1}{2}.
\end{equation}
Assumptions (\ref{Prop1}) and (\ref{Prop2}) are fulfilled in many usual cases. For instance, suppose that $\mathbb{E}
\exp(\beta |Z_n(x)|)<\infty$, for some $\beta >0$. Then, Markov's inequality gives us
\begin{equation}\label{markov}
P(|Z_n|\geq\theta)\leq\mathbb{E}\exp(\beta|Z_n(x)|) \exp(-\beta\theta),\quad \theta\geq 0.
\end{equation}
In the same way, if $Z_n (x)$ satisfies the central limit theorem, then property (\ref{Prop2}) holds. This is the case for the random variables $Z_n (x)$ defined in (\ref{5}).

The problems posed here have been   considered by many   authors for specific sequences of positive linear operators, such as Bernstein and Sz\`asz-Mirakyan operators (see, for instance, Cheng \cite{cheng83}, Zeng and Cheng \cite{zengcheng01},  and  Bustamante et al. \cite{bustamantequesada14}). Important subsets of $B(I)$ and $DB(I)$ are, respectively, the set of functions having bounded variation on $I$, denoted by $BV(I)$, and the set of absolutely continuous functions whose Radon-Nikodym derivative has bounded variation on $I$, denoted by $DBV(I)$. There is a huge literature on approximation results for functions belonging to $BV(I)$ and $DBV(I)$ for particular sequences of operators such as Bernstein, Sz\`asz-Mirakyan, Meyer-K$\ddot{\mbox{o}}$nig and Zeller, King type or exponential type operators (cf. Cheng \cite{cheng84}, Bojanic and Cheng \cite{bojanicheng89}, Guo \cite{guo89}, Zeng and Piriou \cite{zengpiriou98}, $\ddot{\mbox{O}}$zarslan et al. \cite{ozarslan10}, and Abel and Gupta \cite{abelgupta20}, among many others).

Usually, the results referring to the sets $BV(I)$ and $DBV(I)$ take on the following form. If $f\in BV(I)$, denote by $V(f,J)$ the total variation of $f$ on the subset $J\subset I$ and define the function
\begin{align}\label{fx}
f_x (y) &=(f(y) - f(x+)) 1_{I\cap (x,\infty)} (y)\nonumber\\
&+ (f(y) - f(x-)) 1_{I\cap (-\infty , x )} (y),\qquad y\in I.
\end{align}
For the Bernstein polynomials, Zeng and Piriou \cite{zengpiriou98} (see also Bustamante  et al. \cite[Th.17]{bustamantequesada14}) showed that
\begin{eqnarray}\label{Zeng}
   \left|B_n (f,x) - \frac{f(x+) + f(x-)}{2}\right|\leq \dfrac{3}{nx(1-x)+1} \sum_{k=0}^{n-1} V\left( f_x ,J(x,k) \right)\nonumber\\
   + \dfrac{2}{\sqrt{nx(1-x)}+1} \left( |f(x+) - f(x-)| + e_n (x) |f(x) - f(x-)|    \right),
\end{eqnarray}
where
\begin{equation*}J(x,k)=\left[ x-\frac{x}{\sqrt{k+1}} , x+\frac{1-x}{\sqrt{k+1}} \right]\mbox{\quad and\quad }e_n (x) = \sum _{k=1}^{n-1} 1_{\lbrace k/n\rbrace} (x). \end{equation*}
For other results concerning functions in $DBV(I)$, see Section 5.

In this paper, we give approximation results for functions in $B(I)$ and $DB(I)$. The corresponding upper bounds are given in terms of a local first modulus of continuity of the functions under consideration (see Section 3, particularly Definition 3.1). When applied to the subsets $BV(I)$ and $DBV(I)$, such bounds are better and simpler to compute than the total variation term on the right-hand side in (\ref{Zeng}). This is one of the main features of this paper. In fact, we show that our upper bounds are best possible, in the sense that we can construct a sequence of positive linear operators attaining them when acting on large sets of symmetric functions (see Section 4). The main results are given in Section 3. Such results are illustrated in the case of the Bernstein polynomials in Section 5, including a comparative discussion with other  known results in the literature.

We finally point out that the approach given in this paper could be generalized at the price of introducing more involved notations. Specifically, the set $B(I)$ could be enlarged to include unbounded functions, whereas assumptions (\ref{Prop1}) and (\ref{Prop2}) could be weakened (see, for instance, Theorem 4 in Section 3 regarding assumption (\ref{Prop2})).

\section{Preliminary results}\label{preliminaryresults}

Let $\mathcal{L}$ be the set of nondecreasing functions $\psi:[0,\infty ] \to [ 0 , \infty )$ such that
\begin{equation}\label{12}
\psi (0) =\psi(0+) = 0 ,\qquad 0< \psi(\infty) := \lim_{\theta \to \infty} \psi (\theta) < \infty.
\end{equation}
Given $\psi \in \mathcal{L}$, we define
\begin{equation}\label{13}
\Psi (y) = \int_0 ^y \psi (\theta)\, d\theta,\quad y\geq 0.
\end{equation}
On the other hand, let $X_{\beta}$ be a random variable having the exponential distribution with parameter $\beta$, that is,
\begin{equation}\label{14}
F_{\beta}(y):=P(X_{\beta}\leq y)= \int_0^y \beta\, \exp(-\beta \theta)\, d\theta = 1-\exp(-\beta y),\quad y\geq 0.
\end{equation}
A crucial quantity in this paper is the following
\begin{equation}\label{15}
K_n (\beta) := \mathbb{E}\, \psi \left( \frac{X_\beta}{\sqrt{n}}\right),\quad n\in\mathbb{N},\ \psi \in \mathcal{L}.
\end{equation}
By (\ref{12}) and the dominated convergence theorem, we see that $K_n(\beta) = o(1)$, as $n\to \infty$. Upper and lower bounds for $K_n (\beta)$ of the same order of magnitude are given in the following result.
\begin{lemma}\label{lema1}
 Let  $\psi\in\mathcal{L} $  and let $\Psi$ be as in (\ref{13}). For any $m\in\mathbb{N}$, with $m\geq 2$, we have
\begin{equation}\label{16}
\frac{1}{m}\sum_{k=1}^{m-1}\psi\left(\frac{1}{\beta\sqrt{n}}\log\frac{m}{k}\right)\leq \mathbb{E}\,\psi\left(\frac{X_{\beta}}{\sqrt{n}}\right)
\leq \frac{1}{m}\sum_{k=0}^{m-1}\psi\left(\frac{1}{\beta\sqrt{n}}\log\frac{m}{k}\right),
\end{equation}
where it is understood that $\log\infty=\infty$. In addition,
\begin{equation}\label{estima2}
 \mathbb{E}\,\Psi\left(\frac{X_{\beta}}{\sqrt{n}}\right) = \frac{1}{\beta \sqrt{n}}\mathbb{E}\, \psi\left(\frac{X_\beta}{\sqrt{n}}\right).
\end{equation}
\end{lemma}

{\bf Proof.} Let $0=y_0<y_1<\ldots<y_{m-1}<y_m = \infty$  be a finite sequence such that $F_{\beta}(y_j)=j/m$ or, equivalently, $P(y_{j-1}\leq X_{\beta}\leq y_j)=1/m,\, j=1,\ldots , m.$ Note that
\begin{equation*}
y_j = \frac{1}{\beta} \log \frac{m}{m-j},\quad j=1,\ldots, m-1.
\end{equation*}
Since $\psi$ is nondecreasing, we have from (\ref{15})
\begin{equation*}
\mathbb{E}\,\psi\left(\frac{X_{\beta}}{\sqrt{n}}\right)= \sum_{j=1}^m \int_{y_{j-1}}^{y_j} \psi \left( \frac{\theta}{\sqrt{n}} \right) \beta\, \exp(-\beta \theta)\, d\theta \leq\frac{1}{m}  \sum_{j=1}^m \psi    \left(   \frac{y_j}{\sqrt{n}}   \right),
\end{equation*}
which shows the upper bound in (\ref{16}). The lower bound is proved in a similar way.

On the other hand, we have from (\ref{13}), (\ref{14}), and Fubini's theorem
\begin{align*}
\mathbb{E}\, \Psi \left( \frac{X_\beta}{\sqrt{n}}\right)& = \mathbb{E} \int_0^{X_\beta / \sqrt{n}} \psi (\theta)\, d\theta =\mathbb{E} \int_0^{\infty} \psi (\theta) 1_{\lbrace \theta \leq X_\beta / \sqrt{n} \rbrace} \, d\theta
\\
   &= \int_0^\infty \psi (\theta) P(X_\beta \geq \theta \sqrt{n}) d\theta \\
   &= \frac{1}{\beta} \int_0 ^\infty\psi (\theta) \beta \exp(-\beta  \theta\sqrt{n})d\theta = \frac{1}{\beta \sqrt{n}} \mathbb{E} \psi  \left( \frac{X_\beta}{\sqrt{n}}\right),
\end{align*}
thanks to the change $u=\theta\sqrt{n}$. This shows (\ref{estima2}) and completes the proof.\newline \phantom{.} \ \hfill $\square$

Since $K_n (\beta) = o(1)$, as $n\to\infty$, we see from (\ref{estima2}) that
\begin{equation}
\mathbb{E} \Psi \left( \frac{X_\beta}{\sqrt{n}}\right) = \frac{1}{\sqrt{n}} o(1),\quad \text{as }n\to\infty.
\end{equation}
Let $X$ be a nonnegative random variable whose tail probabilities satisfy
\begin{equation}\label{tail}
P(X\geq y) \leq C(\beta)\, \exp(-\beta y),\quad y\geq 0,
\end{equation}
for some positive constant $C(\beta).$
\begin{lemma}\label{2lema}
Let $\psi \in \mathcal{L}$ and let $X$ be as in (\ref{tail}). Then,
\begin{equation*}
\mathbb{E}\psi\left(\frac{X}{\sqrt{n}}\right)\leq C(\beta) \mathbb{E}\, \psi \left(\frac{X_\beta}{\sqrt{n}}\right)
\end{equation*}and
\begin{equation*}
\mathbb{E}\Psi\left(\frac{X}{\sqrt{n}}\right)\leq\frac{C(\beta)}{\beta  \sqrt{n}} \mathbb{E}\, \psi \left(\frac{X_\beta}{\sqrt{n}}\right).\end{equation*}
\end{lemma}
{\bf Proof.} Define
\begin{equation*}\widetilde{\psi} (y) := \lim_{\theta \downarrow y} \psi (y),\quad y\geq 0.
\end{equation*}
The function $\widetilde{\psi} \in \mathcal{L}$ is right-continuous and satisfies
\begin{equation}
\psi (y) \leq \widetilde{\psi}(y),\quad y\geq 0,\qquad \widetilde{\psi} ( \infty) = \psi(\infty).
\end{equation}
Let $W$ be a random variable whose distribution function is given by
\begin{equation}\label{21}
F(y) = P(W\leq y) = \frac{\widetilde{\psi}(y)}{\psi (\infty)},\quad y\geq 0.
\end{equation}
Without loss of generality, we can assume that $W$ is independent of $X$ and $X_\beta$. From (\ref{tail})-(\ref{21}), we infer that
\begin{equation*}\mathbb{E}  \psi \left(\frac{X}{\sqrt{n}}\right) \leq \mathbb{E}\, \widetilde{\psi} \left(\frac{X}{\sqrt{n}}\right)\end{equation*}
\begin{equation}\label{22}
= \psi (\infty) P(W\leq X/\sqrt{n}) \leq \psi (\infty) C(\beta) \mathbb{E}\, \exp(-\beta W\sqrt{n}).
\end{equation}
We repeat the same procedure replacing $X$ by $X_\beta$. Since $\psi$ and $\widetilde{\psi}$ differ in a denumerable set at most and $P(X_\beta \geq y) =\exp(-\beta y),\ y\geq0,$ we get
\begin{equation*} \mathbb{E}  \psi \left(\frac{X_\beta}{\sqrt{n}}\right) = \mathbb{E}  \widetilde{\psi} \left(\frac{X_\beta}{\sqrt{n}}\right) = \psi (\infty) P(W\leq X_\beta/\sqrt{n}) \leq \psi (\infty)   \mathbb{E}\exp(-\beta  W\sqrt{n}). \end{equation*}
This and (\ref{22}) show the first inequality in Lemma \ref{2lema}. The proof of the second one is analogous to that of (\ref{estima2}) and therefore we omit it. \hfill $\square$

\section{Main results}\label{mairesults}
The upper bounds in the approximation results in this paper will be given in terms of the local first modulus of continuity defined as follows.
\begin{definition}\label{definition31}
\emph{Let $x\in\overset{\circ}{I}$ and $h\in B(I)$. The first modulus of continuity of $h$ at $x$ is defined as}
\begin{equation}\label{22*}
\omega_x (h , \theta) = \sup \lbrace |h(x+u)-h(x)|:\, x+u\in I,\, |u|\leq \theta \rbrace , \quad \theta \geq 0.
\end{equation}
\end{definition}
From now on, we fix $x\in\overset{\circ}{I}$. For this reason, we simply write $\omega_x(h,\cdot)=\omega(h,\cdot)$. Also, if $f\in B(I)$, its associated function $f_x$, as defined in (\ref{fx}), is simply denoted by $g=f_x$. Observe that $g$ is continuous at $x$ and $g(x)=0$. It therefore follows from (\ref{22*}) that
\begin{equation}\label{23}
\omega (g , \theta) = \sup \lbrace |g(x+u)|:\, x+u\in I,\, |u|\leq \theta \rbrace , \quad \theta \geq 0.
\end{equation}
Observe that $\omega(g,\cdot)\in\mathcal{L}$. By distinguishing the cases $y<x$, $y=x$, and $y>x$, it can be checked that
\begin{eqnarray}\label{24}
   f(y) - g(y) -  \frac{1}{2} \left( f(x+)  + f(x-)\right) =\left( f(x) - f(x-)\right) 1_{\lbrace x \rbrace} (y)\nonumber\\
      +\left( f(x+) - f(x-)\right) \left(  1_{I\cap (x,\infty)}(y) - \frac{1}{2}  \right) ,\qquad y\in I.
\end{eqnarray}
We give the following result for the operators $L_n$ considered in (\ref{Ln}).

\begin{theorem}\label{teorema3}
Let $f\in B(I)$ and let $Z_n (x)$ be as in  (\ref{Ln}).  Then,
\begin{eqnarray}\label{24*}
\left|L_n (f,x) - \frac{f(x+) + f(x-)}{2}\right|\leq  |f(x) - f(x-)| P(Z_n (x) =0)\nonumber\\
 + |f(x+) - f(x-)| \left| P(Z_n (x)>0) - \frac{1}{2}\right|+ C_n(\beta)\mathbb{E} \omega \left(g , \frac{X_\beta}{\sqrt{n}} \right),
\end{eqnarray}
where $C_n (\beta)$ is defined in (\ref{Prop1}).
\end{theorem}

{\bf Proof.} Replacing $y$ by the random variable $x+Z_n(x) / \sqrt{n}$ in (\ref{24}) and then taking expectations, we obtain
\begin{eqnarray}\label{25}
 L_n (f,x) - \dfrac{f(x+) + f(x-)}{2}=(f(x) - f(x-)) P(Z_n (x) =0) \nonumber\\
      + (f(x+) - f(x-)) \left(P(Z_n (x)>0) - \frac{1}{2}\right) + L_n (g,x).
\end{eqnarray}
 On the other hand, we have from (\ref{Prop1}) and Lemma \ref{2lema} applied to the function $\omega(g,\cdot)\in\mathcal{L}$, as defined in (\ref{23}),
\begin{equation*}
|L_n (g,x)| = \left|\mathbb{E} g \left(x + \frac{Z_n (x)}{\sqrt{n}}   \right)\right| \leq \mathbb{E} \omega \left( g , \frac{|Z_n (x)|}{\sqrt{n}} \right)\leq C_n (\beta) \mathbb{E} \omega \left(g , \frac{X_{\beta}}{\sqrt{n}} \right),
\end{equation*}
which, in conjunction with (\ref{25}), shows the result.
\hfill$\square$

Recently, Mowzer \cite{taqrik24} has considered a sequence of Bernstein-Durrmeyer type operators for which
\begin{equation}\label{26}
P(Z_n(x)>0)\to p(x),\quad P(Z_n(x)<0)\to q(x),\quad n\to\infty,
\end{equation}
where $p(x)+q(x)=1,\ p(x)\not=1/2$. Under condition (\ref{26}), Theorem \ref{teorema3} can be reformulated as follows. We replace identity (\ref{24}) by
\begin{align*}
&f(y)-g(y)-(f(x+)p(x)+f(x-)q(x))\\=&(f(x+)-f(x))(1_{I\cap (x,\infty)}(y)-p(x))\\
+&(f(x-)-f(x))(1_{I\cap (-\infty,x)}(y)-q(x)),\quad y\in I.
\end{align*}
Following the same steps as in the proof of Theorem \ref{teorema3}, we can show the following result.

\begin{theorem}\label{teorema4}
Let $f\in B(I)$ and let $Z_n(x)$ be as in (\ref{Ln}). Under assumption (\ref{26}), we have
\begin{equation*}
|L_n(f,x)-(f(x+)p(x)+f(x-)q(x))|\leq |f(x+)-f(x)||P(Z_n(x)>0)-p(x))|
\end{equation*}
\begin{equation*}
+|f(x-)-f(x)||P(Z_n(x)<0)-q(x))|+C_n(\beta)\mathbb{E}\omega\left(g,\frac{X_{\beta}}{\sqrt{n}}\right).
\end{equation*}
\end{theorem}

Let $\phi\in DB(I)$ as in (\ref{phi}) and let $g:=f_x\in B(I)$ be the function associated to the Radon-Nikodym derivative $f\in B(I)$ of $\phi$. By (\ref{24}), we can write
\begin{align}\label{27}
\phi(y)-\phi(x)&=\int_{x}^{y}(f(u)-g(u))\, du+\int_{x}^{y}g(u)\, du\nonumber\\
&=\frac{1}{2}(f(x+)+f(x-))(y-x)+\frac{1}{2}(f(x+)-f(x-))|y-x|\nonumber\\
&+\int_0^{y-x}g(x+u)\, du,\quad y\in I.
\end{align}

For functions in $DB(I)$, we give the following result.
\begin{theorem}\label{teorema5}
  Let $\phi\in DB(I)$ and let $Z_n(x)$ be as in (\ref{Ln}). Then,
\begin{equation*}
\left|L_n(\phi,x)-\phi(x)-\frac{1}{2\sqrt{n}}(f(x+)+f(x-))\mathbb{E}Z_n(x)\right.
\end{equation*}
\begin{equation*}
\left.-\frac{1}{2\sqrt{n}}(f(x+)-f(x-))\mathbb{E}|Z_n(x)|\right|
\leq \frac{C_n(\beta)}{\beta\sqrt{n}}\mathbb{E}\omega\left(g,\frac{X_{\beta}}{\sqrt{n}}\right).
\end{equation*}
\end{theorem}
{\bf Proof.} Replacing $y$ by the random variable $x+Z_n(x)/\sqrt{n}$ in (\ref{27}) and then taking expectations, we get
\begin{equation*}
L_n(\phi,x)-\phi(x) =\frac{1}{2\sqrt{n}}(f(x+)+f(x-))\mathbb{E}Z_n(x)
\end{equation*}
\begin{equation}\label{28}
 +\frac{1}{2\sqrt{n}}(f(x+)-f(x-))\mathbb{E}|Z_n(x)| +\mathbb{E}\int^{Z_n(x)/\sqrt{n}}_{0}g(x+u)\, du.
\end{equation}
From (\ref{23}), we see that
\begin{equation*}\left|\int_0^yg(x+u)\, du\right|\leq \int_0^{|y|}\omega(g,u)\, du,\ y\in\mathbb{R}.
\end{equation*}
Applying Lemma (\ref{13}) with $\psi=\omega(g,\cdot)$, this implies that
\begin{equation*}
\left|\mathbb{E}\int_{0}^{Z_n(x)/\sqrt{n}}g(x+u)\, du\right|\leq \mathbb{E}\int_{0}^{|Z_n(x)|/\sqrt{n}}\omega(g,u)\, du=\mathbb{E}\Psi\left(\frac{|Z_n(x)|}{\sqrt{n}}\right)
\end{equation*}
\begin{equation*}
\leq \frac{C_n(\beta)}{\beta\sqrt{n}}\mathbb{E}\omega\left(g,\frac{X_{\beta}}{\sqrt{n}}\right).
\end{equation*}
This and (\ref{28}) show the result.\hfill$\square$.

Observe that Theorems \ref{teorema3} and \ref{teorema5} are meaningful, since
\begin{equation*}
\mathbb{E}\omega\left(g,\frac{X_{\beta}}{\sqrt{n}}\right)=o(1),\quad n\rightarrow\infty,
\end{equation*}
as asserted after definition (\ref{15}). On the other hand, such theorems can be respectively applied to functions in the sets $BV(I)$ and $DBV(I)$ without any modification. Suppose, for instance, that $f\in BV(I)$. It follows from (\ref{23}) that
\begin{equation}\label{29}
  \omega(g,\theta)\leq V(g,[x-\theta,x+\theta]\cap I),\quad \theta\geq 0.
\end{equation}
Therefore, there is no need to express the upper bounds in Theorems \ref{teorema3} and \ref{teorema5} in terms of the total variation of $g$ on appropriate subintervals of $I$, as done in (\ref{Zeng}). More details will be given in Section 5.
\section{The convolution operators $L_n^{\beta}$}\label{section4}
The purpose of this section is to construct a sequence $L_n^{\beta}$ of positive linear operators attaining the upper bounds in Theorems \ref{teorema3} and \ref{teorema5} when acting on certain subsets of symmetric functions.

Let $\beta>0$. Denote by $Y_{\beta}$ a random variable having the symmetric exponential density
\begin{equation*}
\rho_{\beta}(\theta)=\frac{\beta}{2}\, \exp(-\beta|\theta|),\quad \theta \in \mathbb{R}.
\end{equation*}
Note that
\begin{equation}\label{30}
\mathbb{E}h(Y_{\beta})=\frac{\beta}{2}\int_{0}^{\infty}(h(\theta)+h(-\theta))\, \exp(-\beta\theta)\, d\theta=\frac{1}{2}\mathbb{E}(h(X_{\beta})+h(-X_{\beta})),
\end{equation}
for any measurable function \mbox{$h:\mathbb{R}\rightarrow\mathbb{R}$} for which the preceding integral makes sense, $X_{\beta}$ being the random variable defined in (\ref{14}).

We consider the sequence $L_n^{\beta}$ of convolution operators defined as
\begin{equation}\label{31}
L_n^{\beta}(f,x)=\mathbb{E}f\left(x+\frac{Y_{\beta}}{\sqrt{n}}\right),\quad x\in\mathbb{R},\ n\in\mathbb{N},
\end{equation}
acting on measurable functions $f:\mathbb{R}\rightarrow\mathbb{R}$ for which the preceding expectations exist. As follow from (\ref{30}) and (\ref{31}), we can write
\begin{equation}\label{32}
L_n^{\beta}(f,x)=\frac{1}{2}\left(\mathbb{E}f\left(x+\frac{X_{\beta}}{\sqrt{n}}\right)+\mathbb{E}f\left(x-\frac{X_{\beta}}{\sqrt{n}}\right)\right),\quad x\in\mathbb{R},\quad n\in\mathbb{N.}
\end{equation}

Fix $x\in\mathbb{R}$. Denote by $\mathcal{S}\subset B(\mathbb{R})$ the set of functions $f$ which are symmetric around $x$ and nondecreasing in $[x,\infty)$. Note that
\begin{equation}\label{33}
a:=f(x+)=f(x-),\qquad f\in\mathcal{S}.
\end{equation}
Also, it follows from (\ref{32}) that
\begin{equation}\label{34}
L_n^{\beta}(f,x)=\mathbb{E}f\left(x+\frac{X_{\beta}}{\sqrt{n}}\right),\qquad f\in\mathcal{S}.
\end{equation}
Finally, denote by $\widetilde{\mathcal{S}}$ the set of symmetric functions $\widetilde{\phi}$ having the form
\begin{equation}\label{35}
\widetilde{\phi}(y)=\int_0^{|y-x|}f(x+u)\, du,\quad y\in\mathbb{R},\quad f\in\mathcal{S}.
\end{equation}
We state the main result of this section.
\begin{theorem}\label{teorema6}
Let $a$ be as in (\ref{33}). For the operators $L_n^{\beta}$ defined in (\ref{31}), we have

{ (a)} If $f\in\mathcal{S}$, there is equality in Theorem \ref{teorema3}. Specifically,
\begin{equation*}L_n^{\beta}(f,x)-a=\mathbb{E}\omega\left(g,\frac{X_{\beta}}{\sqrt{n}}\right).
\end{equation*}
\quad (b) If $\widetilde{\phi}\in\widetilde{\mathcal{S}}$, there is equality in Theorem \ref{teorema5}. More precisely,
\begin{equation*}L_n^{\beta}(\widetilde{\phi},x)-\frac{a}{\beta\sqrt{n}}=\frac{1}{\beta\sqrt{n}}\mathbb{E}\omega\left(g,\frac{X_{\beta}}{\sqrt{n}}\right).
\end{equation*}

\end{theorem}
{\bf Proof.} {\bf(a)} In the first place, since $Y_{\beta}=Z_n(x)$ is an absolutely continuous and symmetric random variable, we have
\begin{equation*}P(Y_{\beta}=0)=0,\quad P(Y_{\beta}>0)=\frac{1}{2},\quad C_n(\beta)=1.\end{equation*}
By (\ref{fx}) and (\ref{33}), the function $g$ associated to $f$ is given by
\begin{equation*}g(y)=(f(y)-a)1_{\mathbb{R}\setminus\{x\}}(y),\quad y\in\mathbb{R}.\end{equation*}
Since $f\in\mathcal{S}$, this implies, by virtue of (\ref{33}), that
\begin{equation}\label{36}
\omega(g,\theta)=f(x+\theta)-a,\quad \theta\geq 0.
\end{equation}
In turn, this implies that
\begin{equation*}\mathbb{E}\omega\left(g,\frac{X_{\beta}}{\sqrt{n}}\right)=\mathbb{E}f\left(x+\frac{X_{\beta}}{\sqrt{n}}\right)-a=L_n^{\beta}(f,x)-a,\end{equation*}
where the last equality follows from (\ref{34}). This shows part (a).

{\bf (b)} Taking $h(y)=y$ and $h(y)=|y|$ in (\ref{30}), we see that
\begin{equation*}\mathbb{E}Y_{\beta}=0,\quad \mathbb{E}|Y_{\beta}|=\mathbb{E}X_{\beta}=\frac{1}{\beta}.\end{equation*}
Since the function $\widetilde{\phi}$ defined in (\ref{35}) is symmetric around $x$, we have from (\ref{32}) and (\ref{36})
\begin{equation*}L_n^{\beta}(\widetilde{\phi},x)=\mathbb{E}\widetilde{\phi}\left(x+\frac{X_{\beta}}{\sqrt{n}}\right)=\mathbb{E}\int_0^{X_{\beta}/\sqrt{n}}f(x+u)\, du\end{equation*}
\begin{equation*}=\frac{a}{\sqrt{n}}\mathbb{E}X_{\beta}+ \mathbb{E}\int_0^{X_{\beta}/\sqrt{n}}\omega(g,u)\, du=\frac{a}{\beta\sqrt{n}}+
\frac{1}{\beta\sqrt{n}}\mathbb{E}\omega\left(g,\frac{X_{\beta}}{\sqrt{n}}\right),\end{equation*}
where the last equality follows from (\ref{estima2}) for $\psi=\omega(g,\cdot)$. This shows part (b) and completes the proof.\hfill$\square$

\section{Bernstein polynomials}
We illustrate Theorems \ref{teorema3} y \ref{teorema5} in the case of the Bernstein polynomials. To give explicit approximation results, some auxiliary lemmas will be needed. Recall that  Stirling's approximation states that
\begin{equation}\label{37}
{\mbox{$\sqrt{2\pi k}\, \exp\left(\frac{1}{12k+1}\right)\leq \dfrac{k!}{k^k}\, \exp(k)\leq     \sqrt{2\pi k}\, \exp\left(\frac{1}{12k}\right),\ k\in\mathbb{N}.
$}}
\end{equation}
\begin{lemma}\label{lema7}
Let $Z_n(x)$ be as in (\ref{5}). Then,
\begin{equation*}P(Z_n(x)=0)\leq a_n(x):=\frac{1}{\sqrt{2\pi nx(1-x)}}1_{A_n}(x),\, A_n=\left\{\frac{k}{n}:k=1,\ldots,n-1\right\}.\end{equation*}
\end{lemma}
{\bf Proof.} By (\ref{leybinomial}) and (\ref{5}), we have
\begin{equation}\label{38}
P(Z_n(x)=0)=P(S_n(x)=nx)=\sum_{k=1}^{n-1}P(S_n(x)=k)1_{\{k/n\}}(x).
\end{equation}
Differentiating with respect to $x$, it can be checked that
\begin{equation*}
P(S_n(x)=k)\leq P(S_n(k/n)=k),\quad k=1,\ldots,n-1.
\end{equation*}
By (\ref{37}) and the fact that  $1/(12n) -1/(12k+1)-1/(12(n-k)+1)<0$, this implies that
\begin{equation*}P(S_n(x)=k)\leq \frac{1}{\sqrt{2\pi}}\sqrt{\frac{n}{k(n-k)}},\quad k=1,\ldots,n-1.
\end{equation*}
This and (\ref{38}) show the result.\hfill $\square$

Let $Z$ be a random variable having the standard normal density and let $(X_k)_{k\geq 1}$ be a sequence of independent copies of a random variable $X$ such that $\mathbb{E}X=\mu,\ 0<\mathbb{E}(X-\mu)^2=\sigma^2$. and $\mathbb{E}|X-\mu|^3=\gamma<\infty.$ The Berry-Esseen bounds establish that
\begin{equation}\label{39}
\sup_{y\in\mathbb{R}}\left|P\left(\frac{X_1+\ldots+X_n-n\mu}{\sigma\sqrt{n}}\leq y\right)-P(Z\leq y)\right|\leq C\frac{\gamma}{\sigma^3}\frac{1}{\sqrt{n}},
\end{equation}
for some positive constant $C$. According to Essen \cite{esseen56}, $C$ cannot be less than $(3+\sqrt{10})/(\sigma\sqrt{2\pi})=0.4097\cdots$. Shevtsova \cite{shevtsova14} (see also Korolev and Shevtsova \cite{korolev10}) showed that $C\leq 0.4690.$

\begin{lemma}\label{lema8}
Let $Z_n(x)$ be as in (\ref{5}). Then,
\begin{equation*}\left|P(Z_n(x)>0)-\frac{1}{2}\right|\leq b_n(x):=\left|\frac{1}{2}-(1-x)^n\right|1_{\left(0,\frac{1}{n}\right)}(x)\end{equation*}
\begin{equation*}+\left|\frac{1}{2}-x^n\right|1_{\left.\left[\frac{n-1}{n},1\right.\right)}(x)+0.4690\frac{x^2+(1-x)^2}{\sqrt{n x(1-x)}}1_{\left.\left[\frac{1}{n},\frac{n-1}{n}\right.\right)}(x).\end{equation*}
\end{lemma}
{\bf Proof.} If\quad $0<x<1/n$\quad or\quad $(n-1)/n\leq x<1,$ the upper bound $b_n(x)$ is easily checked. Suppose that $1/n\leq x<(n-1)/n$. We apply the Berry-Esseen bounds in (\ref{39}) to the random variable $S_n(x)$ defined in (\ref{2*}). In this regard, note that
\begin{equation*}\sigma^2=\mathbb{E}|S_1(x)-x|^2=x(1-x),
\end{equation*}
and
\begin{equation}\label{40}
\gamma=\mathbb{E}|S_1(x)-x|^3=x(1-x)(x^2+(1-x)^2).
\end{equation}
Since
\begin{equation*}\left|P(Z_n(x)>0)-\frac{1}{2}\right|=\left|P\left(\frac{S_n(x)-nx}{\sqrt{nx(1-x)}}\leq 0\right)-P(Z\leq 0)\right|,\end{equation*}
the result follows from (\ref{39}), (\ref{40}), and the comments following (\ref{39}).\hfill $\square$

To estimate the constant $C_n(\beta)$ in (\ref{Prop1}) referring to the Bernstein polynomials, we give the following result.
\begin{lemma}\label{lema9}
Let $Z_n(x)$ be as in (\ref{5}). For any $\beta>0$, we have
\begin{equation*}\mathbb{E}\exp(\beta|Z_n(x)|)\leq 2\, \exp\left(\frac{x(1-x)}{2}\left(\beta^2+\frac{x^2+(1-x)^2}{3\sqrt{n}}\beta^3\, \exp\left(\frac{\beta}{\sqrt{n}}\right)\right)\right).\end{equation*}
\end{lemma}
{\bf Proof.} The result will follow as soon as we show that
\begin{equation*}
\mathbb{E}\exp(\beta|S_n(x)-nx|)
\end{equation*}\begin{equation}\label{41}
\leq 2\, \exp\left(\frac{nx(1-x)}{2}\left(\beta^2+\frac{x^2+(1-x)^2}{3}\beta^3\, \exp(\beta)\right)\right),
\end{equation}
since, by (\ref{5}), it will suffice to replace $\beta$ by $\beta/\sqrt{n}$ in (\ref{41}). For $k=3,4,\ldots,$ we have
\begin{equation*}|\mathbb{E}(S_1(x)-x)^k|=x(1-x)|(1-x)^{k-1}+(-1)^kx^{k-1}|\leq x(1-x)(x^2+(1-x)^{2}).\end{equation*}
Therefore,
\begin{equation*}\mathbb{E}\exp(\beta(S_1(x)-x))\leq 1+\frac{x(1-x)}{2}\beta^2+\sum_{k=3}^{\infty}\frac{\beta^k}{k!}|\mathbb{E}(S_1(x)-x)^k|\end{equation*}
\begin{equation*}\leq 1+x(1-x)\left(\frac{\beta^2}{2}+(x^2+(1-x)^2)\sum_{k=3}^{\infty}\frac{\beta^k}{k!}\right)\end{equation*}
\begin{equation*}
\leq \exp\left(x(1-x)\left(\frac{\beta^2}{2}+(x^2+(1-x)^2)\left( \exp(\beta)-1-\beta-\frac{\beta^2}{2}\right)\right)\right).\end{equation*}
Since $\exp(\beta)-1-\beta-\beta^2/2\leq \beta^3\, \exp(\beta)/6$, this implies that
$$\mathbb{E}\exp(\beta(S_n(x)-nx))=\left(\mathbb{E} \exp(\beta(S_1(x)-x))\right)^n$$
\begin{equation}\label{42}
  \leq \exp\left(\frac{nx(1-x)}{2}\left(\beta^2+\frac{x^2+(1-x)^2}{3}\beta^3\, \exp(\beta)\right)\right).
\end{equation}
The random variables $S_n(1-x)$ and $n-S_n(x)$ have the same law. Hence,
\begin{equation*}
\mathbb{E}\exp(\beta|S_n(x)-nx|)=\mathbb{E}\exp(\beta(S_n(x)-nx))1_{\{S_n(x)\geq nx\}}\end{equation*}
\begin{equation*}+\mathbb{E}\exp(-\beta(S_n(x)-nx))1_{\{S_n(x)< nx\}}
\end{equation*}
\begin{equation*}
\leq \mathbb{E}\exp(\beta(S_n(x)-nx))+\mathbb{E}\exp(\beta(S_n(1-x)-n(1-x))).
\end{equation*}
This and (\ref{42}) shows claim  (\ref{41}) and completes the proof.\hfill $\square$

We are in a position to apply Theorem \ref{teorema3} and \ref{teorema5} to the case of the Bernstein polynomials.
\begin{corollary}\label{corolario10}
Let $f\in B([0,1])$ and $x\in (0,1)$. Then,
\begin{equation*}
\left|B_n(f,x)-\frac{f(x+)+f(x-)}{2}\right|\leq a_n(x)|f(x)-f(x-)|+b_n(x)|f(x+)-f(x-)|
\end{equation*}
\begin{equation*}
+2 \exp\left(\frac{x(1-x)}{2}\left(1+\frac{x^2+(1-x)^2}{3\sqrt{n}}\, \exp\left(\frac{1}{\sqrt{n}}\right)\right)\right)
\end{equation*}
\begin{equation}\label{43}
\times \frac{1}{n}\sum_{k=0}^{n-1}\omega\left(g,\frac{1}{\sqrt{n}}\log \frac{n}{k}\right),
\end{equation}
where $a_n(x)$ and $b_n(x)$ are defined in Lemmas \ref{lema7} and \ref{lema8}, respectively.
\end{corollary}
{\bf Proof.} Starting from Theorem \ref{teorema3}, use Lemmas \ref{lema7} and \ref{lema8} to bound the terms multiplying the quantities $|f(x)-f(x-)|$ and $|f(x+)-f(x-)|$, respectively. To bound the last term on the right-hand side in (\ref{24*}), use (\ref{markov}) and Lemma \ref{lema9} with $\beta=1$
to estimate $C_n(\beta)$, and inequality (\ref{16}) with $\psi=\omega(g,\cdot)$, $m=n$, and $\beta=1$ to bound the term
\begin{equation*}
\mathbb{E}\omega\left(g,\frac{X_{\beta}}{\sqrt{n}}\right).
\end{equation*}
This completes the proof.\hfill $\square$

The comparison between Corollary \ref{corolario10} and formula (\ref{Zeng}) reveals the following. The terms $a_n(x)$ and $b_n(x)$ in (\ref{43}) are smaller than the corresponding ones in expression (\ref{Zeng}). However, the main difference lies in the last term on the right-hand side in (\ref{43}). On the one hand, it follows from (\ref{29}) that
\begin{eqnarray}\label{44}
&\displaystyle{\sum_{k=0}^{n-1}\omega\left(g,\frac{1}{\sqrt{n}}\log\frac{n}{k}\right)}&\nonumber\\
\leq & \displaystyle{\sum_{k=0}^{n-1}V\left(g,\left[x-\frac{1}{\sqrt{n}}\log\frac{n}{k},x+\frac{1}{\sqrt{n}}\log\frac{n}{k}\right]\cap[0,1]\right)},&
\end{eqnarray}
and the term on the left is easier to compute than term on the right in (\ref{44}). On the other hand, excepting the case $k=0$, the lengths of the intervals in (\ref{44}) are asymptotically much shorter than those in (\ref{Zeng}). For instance, the lengths of the $(n-1)th$ intervals in (\ref{Zeng}) and (\ref{44}) are, respectively,
\begin{equation*}\frac{1}{\sqrt{n}}\mbox{\quad and\quad  }\frac{2}{\sqrt{n}}\log\frac{n}{n-1}\sim \frac{2}{n\sqrt{n}},\quad n\rightarrow\infty.
\end{equation*}

Denote by $\lfloor y\rfloor$ and $\lceil y\rceil$ the floor and the ceiling of $y\in\mathbb{R}$, respectively.
\begin{corollary}\label{corolario11}
Let $\phi\in DB([0,1])$ and $x\in(0,1)$. Then,
\begin{equation*}
|B_n(\phi,x)-\phi(x)-x(1-x)P(S_{n-1}(x)=\lfloor nx\rfloor)(f(x+)-f(x-))|
\end{equation*}
\begin{equation*}
\leq 2 \exp\left(\frac{x(1-x)}{2}\left(1+\frac{x^2+(1-x)^2}{3\sqrt{n}}\, \exp\left(\frac{1}{\sqrt{n}}\right)\right)\right)
\end{equation*}
\begin{equation}\label{45}
\times \frac{1}{n}\sum_{k=0}^{\lceil \sqrt{n}\rceil-1}\omega\left(g,\frac{1}{\sqrt{n}}\log \frac{n}{k}\right).
\end{equation}
\end{corollary}
{\bf Proof.} It is known (cf. Johnson and Kotz  \cite[p. 52]{johnson69}) that
\begin{equation}\label{46}
\mathbb{E}\left|\frac{S_n(x)-nx}{n}\right|=2x(1-x)P(S_{n-1}(x)=\lfloor nx\rfloor).
\end{equation}
Hence, starting from Theorem \ref{teorema5} and recalling (\ref{5}), we see that

\begin{equation*}
\mathbb{E}Z_n(x)=0,\quad \frac{1}{2\sqrt{n}}\mathbb{E}|Z_n(x)|=x(1-x)P(S_{n-1}(x)=\lfloor nx \rfloor ).
\end{equation*}
As in the proof of Corollary \ref{corolario10}, the upper bound is (\ref{45}) is obtained by choosing $\beta=1$ in Lemma \ref{lema9}, and $\beta=1$ and $m=\lceil \sqrt{n}\rceil$ in (\ref{16}). This completes the proof.\hfill $\square$

Bojanic and Cheng \cite[Th. 2]{bojanicheng89}, showed the following result. In the setting of Corollary \ref{corolario11}, assume further that $n\geq x(1-x)$. Then,
\begin{equation*}
\left|B_n(\phi,x)-\phi(x)-\frac{1}{\sqrt{n}}\left(\frac{x(1-x)}{2\pi}\right)^{1/2}(f(x+)-f(x-))\right|
\end{equation*}
\begin{equation*}
\leq \frac{M}{2n\sqrt{x(1-x)}}|f(x+)-f(x-)|
\end{equation*}
\begin{equation}\label{47}
+\frac{2}{n}\sum_{k=0}^{\lfloor \sqrt{n}-1\rfloor}V\left(g,\left[x-\frac{x}{\sqrt{k+1}},x+\frac{1-x}{\sqrt{k+1}}\right]\right),
\end{equation}
where $M$ is an explicit constant (cf.\cite[Th.1 ]{bojanicheng89}). We point out that there is a misprint in \cite[Th 2]{bojanicheng89}, concerning the main term of the approximation.

Observe that there is no restriction on $n$ in Corollary \ref{corolario11}. On the other hand, the main term of the approximation in (\ref{45}) and (\ref{47}) is essentially the same. Indeed, we have from (\ref{46}) and the central limit theorem for $S_n(x)$
\begin{equation*}
x(1-x)P(S_{n-1}(x)=\lfloor nx\rfloor)=\frac{\sqrt{x(1-x)}}{2\sqrt{n}}\mathbb{E}\left|\frac{S_n(x)-nx}{\sqrt{nx(1-x)}}\right|
\end{equation*}
\begin{equation*}
\sim \frac{\sqrt{x(1-x)}}{2\sqrt{n}}\mathbb{E}|Z|=\frac{1}{\sqrt{n}}\left(\frac{x(1-x)}{2\pi}\right)^{1/2},\quad n\rightarrow\infty,
\end{equation*}
where $Z$ is a standard normal random variable.
Finally, the asymptotic value of the constant in the upper bound in (\ref{45}) is $2\,\exp(x(1-x)/2),$ which is slightly worse than $2$. However, the total variation term in (\ref{47}) is worse than the upper bound in (\ref{45}), as explained in the comments following Corollary \ref{corolario10}.
\section*{Acknowledgments}
The first author is supported by Research Project DGA (E48\_23R). The second and third authors are supported by Junta de Andaluc\'ia Research Group  FQM-0178.

\end{document}